\documentclass[12pt,a4paper,leqno]{article}

\usepackage{amssymb,exscale}
\usepackage[centertags]{amsmath}
\usepackage{amsthm}

\numberwithin{equation}{section}

\swapnumbers
\theoremstyle{definition}
\newtheorem{theorem}[equation]{Theorem}
\newtheorem{lemma}[equation]{Lemma}

\newtheorem{definition}[equation]{Definition}

\renewcommand{\phi}{\varphi}

\renewcommand{\(}{\bigl(}
\renewcommand{\)}{\bigr)\vphantom{)}}

\newcommand{\spec}{\operatorname{spec}}

\newcommand{\Proj}{\operatorname{Pr}}

\newcommand{\bN}{\mathbf N}
\newcommand{\Om}{\Omega}
\newcommand{\om}{\omega}

\newcommand{\s}{\sigma}

\newcommand{\F}{\mathcal F}

\newcommand{\A}{\mathcal A}
\newcommand{\la}{\lambda}

\newcommand{\Fstable}{\F_{\text{stable}}}
\newcommand{\Ex}{\mathbb E}

\newcommand{\R}{\mathbb R}

\newcommand{\Z}{\mathbb Z}

\newcommand{\cE}[2]{\mathbb{E}\,\(\,#1\,\big|\,#2\,\)\,}
\newcommand{\CE}[2]{\mathbb{E}\,\bigg(#1\,\bigg|\,#2\bigg)\,}

\newcommand{\sif}{$\sigma$-field}

\newcommand{\One}{\mathbf1}
\newcommand{\imply}{\;\;\;\Longrightarrow\;\;\;}

\newcommand{\imp}{$ \Longrightarrow $ }
\newcommand{\impl}{\;\Longrightarrow\;}
\renewcommand{\equiv}{\;\;\;\Longleftrightarrow\;\;\;}

\def\emailwww#1#2{\par\qquad {\tt #1}\par\qquad {\tt #2}\medskip}

\begin{document}

\title{Noise sensitivity on continuous products:\\
an answer to an old question of J.~Feldman}
\author{Boris Tsirelson}
\date{}

\maketitle

\begin{abstract}
A relation between $\s$-additivity and linearizability, conjectured by
Jacob Feldman in 1971 for continuous products of probability spaces,
is established by relating both notions to a recent idea of noise
stability/sensitivity.
\end{abstract}

\section*{Introduction}
A discrete-time random process with independent values is just a
sequence of independent random variables, described by the product of
a sequence of probability spaces. What could be its continuous-time
counterpart? Non-equivalent approaches were proposed
\cite{Fe,unitary,TV}, the whole picture being still unclear.

Independent \sif s are a more convenient language than products of
probability spaces. Each approach deals with a family $
(\F_A)_{A\in\A} $ of sub-\sif s $ \F_A \subset \F $ on a probability
space $ (\Om,\F,P) $, indexed by subsets $ A \subset T $ of some
``base'' set $ T $, belonging to an algebra\footnote{%
That is, $ A \in \A \impl T \setminus A \in \A $ and $ A,B \in \A
\impl A \cup B \in \A $.}
$ \A $; the family satisfies
\begin{equation}\label{0.1}
\F_{A \uplus B} = \F_A \otimes \F_B \, .
\end{equation}
That is, if $ A,B \in \A $ and $ A \cap B = \emptyset $, then $ \F_A $
and $ \F_B $ are independent\footnote{%
It means that $ P (E\cap F) = P(E) P(F) $ for all $ E \in \F_A $, $ F
\in \F_B $.}
and, taken together, they generate $ \F_{A\cup B} $. Approaches differ
in additional conditions on $ \A $ and $ (\F_A) $. Most restrictive
approaches admit (generalized versions of) classical results such as
Levy-Khintchine formula and Levy-Ito theorem. Less restrictive
approaches (at least, some of them) are not at all pathologic, they
arise from quite natural finite models whose scaling limits go beyond
the classical theory \cite{scaling,TV,W}.

The approach used by Feldman in 1971 \cite{Fe} requires $ \A $ to be
the Borel \sif\ of a standard Borel space, and $ (\F_A) $ to be
$\s$-additive in the sense that\footnote{%
In other words, if $ A_1 \subset A_2 \subset \dots $ and $ A = A_1
\cup A_2 \cup \dots $ then $ \F_A $ is generated by $ \F_{A_1} \cup
\F_{A_2} \cup \dots $}
\begin{equation}\label{0.2}
A_n \uparrow A \imply F_{A_n} \uparrow \F_A \, .
\end{equation}
The classical theory holds \cite{Fe} for every \emph{decomposable
process,} defined as a family $ (X_A)_{A\in\A} $ of random
variables\footnote{%
A random variable is treated as an equivalence class of measurable
functions on $ \Om $.}
$ X_A $ such that
\begin{equation}
\begin{gathered}
X_A \text{ is $\F_A$-measurable} \, , \\
X_{A\uplus B} = X_A + X_B \, , \\
A_n \uparrow A \imply X_{A_n} \to X_A \, .
\end{gathered}\end{equation}
The problem formulated by Feldman \cite[Problem 1.9]{Fe}: (a) Does
every $ (\F_A) $ possess a nontrivial decomposable process? More
strongly: (b) Is every $ (\F_A) $ \emph{linearizable,} that is,
generated by its decomposable processes?\footnote{%
Feldman treats a decomposable process more generally (it is defined
on some ideal, not the whole $ \A $). We do not need it, since $
(\F_A) $ is generated by decomposable processes defined on the whole $
\A $.}
Both questions are answered below in the positive. To this end, a
concept of noise stability/sensitivity \cite{BKS} will be adapted to
the continuous case.

Feldman's framework is quite restrictive in demanding $ \A $ to be a
\sif. Recent examples \cite{scaling,TV,W} provide $ \F_A $ only for
elementary sets $ A $, that is, finite unions of intervals. (Extending
$ (\F_A) $ to more general $ A $ is often impossible, as will be
seen.)  Restricting ourselves to intervals with rational endpoints we
get a \emph{countable} algebra $ \A $ of sets, which is a convenient
framework, used in Sections 2, 3.

\section{Elementary case}
In this section the algebra $ \A $ is assumed to be finite. Thus, $ \A
$ corresponds to a finite partition $ T = a_1 \uplus \dots \uplus a_m
$, and
\[
\F_T = \F_{a_1} \otimes \dots \otimes \F_{a_m} \, .
\]
Each $ A \in \A $ is of the form $ A = a_{k_1} \uplus \dots \uplus
a_{k_n} $, and \eqref{0.1} means simply $ \F_A = \F_{a_{k_1}} \otimes
\dots \otimes \F_{a_{k_n}} $. Ascribing to $ A $ the
probability\footnote{%
That is \emph{not} the probability $ P $ appearing in $ (\Om,\F,P) $.}
\begin{equation}\label{1.1}
\mu_p (A) = \mu_p ( a_{k_1} \uplus \dots \uplus a_{k_n} ) = p^n
(1-p)^{m-n}
\end{equation}
we get Bernoulli measure $ \mu_p $ on $ \A $; $ p \in [0,1] $ is its
parameter. Note that $ \mu_p $ is not a measure on $ (T,\A) $, it is
rather a measure on $ (\A, 2^\A) $; in other words, $ \A $ is treated
here as just a set (not an algebra), equipped with the \sif\ $ 2^\A $
of all its subsets.\footnote{%
Thus, $ \mu_p(A) $ should be written rigorously as $ \mu_p \( \{A\} \)
$.}

Imagine that $ A,B \in \A $ are chosen at random, independently,
according to $ \mu_{p_1} $ and $ \mu_{p_2} $ respectively; then $ A
\cap B $ is a random set distributed $ \mu_{p_1 p_2} $. In other
words,
\begin{equation}\label{1.2}
\mu_{p_1} * \mu_{p_2} = \mu_{p_1 p_2} \, ;
\end{equation}
here the convolution ($ * $) of measures on $ \A $ (that is, on $ (
\A, 2^\A ) $) is taken with respect to the semigroup operation of
intersection, $ \A \times \A \ni (A,B) \mapsto A \cap B \in \A $. The
corresponding continuous-time Markov process on $ \A $ (its time $ t $
is related to $ p $ by $ p = e^{-t} $) is easy to describe; initially
(at $ t = 0 $) the random set is the whole $ T $; during an
infinitesimal time interval $ (t,t+dt) $ each $ a_k $ is excluded from
the random set with probability $ dt $; choices are independent for $
k = 1,\dots,m $; if $ a_k $ was excluded before, nothing happens.

A conditional expectation operator corresponds to every $ A \in \A $,
\begin{equation}\label{1.3}
E_A : L_2 (\F_T) \to L_2 (\F_T) \, , \qquad E_A (X) = \cE{ X }{ \F_A }
\, , \quad E_A = \operatorname{Pr}_{L_2(\F_A)} \, ,
\end{equation}
just the orthogonal projection onto\footnote{%
$ L_2 (\F_A) $ means $ L_2 (\Om,\F_A,P) $.}
$ L_2 (\F_A) \subset L_2 (\F_T) $. Note that $ E_T = \One $ (since the
operators act on $ L_2 (\F_T) $, not the whole $ L_2 (\F) $), and
\begin{equation}\label{1.4}
E_A E_B = E_{A\cap B} \, ,
\end{equation}
however, $ E_{A \uplus B} $ is not $ E_A + E_B $; $ (E_A)_{A\in\A} $
is not a projection measure on $ (T, \A) $.
In order to get a joint diagonalization of the commuting operators $
E_A $, introduce for every $ A \in \A $ a space $ H_A $ consisting of
all $ X \in L_2 (\F_A) $ that are orthogonal to $ L_2 (\F_B) $ for all
$ B \subset A $, $ B \ne A $. We have
\begin{equation}\label{1.5}
\begin{gathered}
L_2 (\F_T) = \bigoplus_{A\in\A} H_A \, ; \\
X = \sum_{A\in\A} X_A \, , \;\; X_A \in H_A \imply E_A X = \sum_{B
  \subset A} X_B \, .
\end{gathered}
\end{equation}
Note that $ H_\emptyset = L_2 (\F_\emptyset) $ is the one-dimensional
space of constants, and
\begin{equation}\label{1a}
H_{A\uplus B} = H_A \otimes H_B
\end{equation}
in the sense that, for any two disjoint $ A,B \in \A $, random
variables of the form $ XY $ for $ X \in H_A $, $ Y \in H_B $ (belong
to and) span $ H_{A\cup B} $. In other words,
\begin{equation}\label{1b}
H_{a_{k_1} \uplus \dots \uplus a_{k_n}} = H_{a_{k_1}} \otimes \dots
\otimes H_{a_{k_n}} \, ;
\end{equation}
a proof for $ H_{a_1\cup a_2} $ (general case being similar) consists
in choosing orthogonal bases $ (X_i)_i $ in $ L_2 (\F_{a_1}) $, $ X_0
= \One $, and $ (Y_j)_j $ in $ L_2 (\F_{a_2}) $, $ Y_0 = \One $, and
considering the basis $ (X_i Y_j)_{i,j} $ in $ L_2 (\F_{a_1}) \otimes
L_2 (\F_{a_2}) = L_2 (\F_{a_1} \otimes \F_{a_2}) = L_2 ( \F_{a_1\cup
a_2} ) $. So,
\begin{equation}\label{1b1}
\begin{split}
L_2 (\F_T) &= L_2 (\F_{a_1}) \otimes \dots \otimes L_2 (\F_{a_m}) = \\
&= \(
  H_\emptyset \oplus H_{a_1} \) \otimes \dots \otimes \( H_\emptyset
  \oplus H_{a_m} \) = \bigoplus_{A\in\A} H_A \, , \\
\Proj_{H_A} &= \bigg( \bigotimes_{a\subset A} \Proj_{H_a} \bigg)
  \otimes \bigg( \bigotimes_{a\subset T\setminus A} \( \One -
  \Proj_{H_a} \) \bigg) \, .
\end{split}\end{equation}

Combining the conditional expectations with the convolution semigroup,
we get an operator semigroup
\begin{gather}
U_t : L_2 (\F_T) \to L_2 (\F_T) \, , \notag \\
U_t = \int E_A \, d\mu_p (A) = \sum_{A\in\A} \mu_p(A) E_A \quad
\text{where } p = e^{-t} \, , \label{1.6} \\
U_s U_t = U_{s+t} \, , \qquad U_0 = \One \, . \notag
\end{gather}
In the language of tensor products,
\begin{multline}\label{1c}
U_t = \( \underbrace{ \One \oplus e^{-t} \cdot \One }_{\text{on }
  H_\emptyset \oplus H_{a_1}} \) \otimes \dots \otimes \(
  \underbrace{ \One \oplus e^{-t} \cdot \One }_{\text{on }
  H_\emptyset \oplus H_{a_m}} \) = \\
= \sum_{A\in\A} \bigg( \bigotimes_{a\subset A} e^{-t} \cdot \One
  \bigg) \otimes \bigg( \bigotimes_{a\subset T\setminus A} \One \bigg)
  \, ,
\end{multline}
and we get eigenspaces
\begin{equation}\label{1d}
\begin{gathered}
H_n = \bigoplus \{ H_A : A = a_{k_1} \uplus \dots \uplus a_{k_n}, \,
  k_1 < \dots < k_n \} \, ; \\
H_0 = H_\emptyset = L_2 (\F_\emptyset) = \text{constants} \, , \quad
  H_m = H_T = H_{a_1} \otimes \dots \otimes H_{a_m} \, , \\
L_2 (\F_T) = H_0 \oplus \dots \oplus H_m \, ; \\
X \in H_n \imply U_t X = e^{-nt} X \, , \\
\spec U_t = \{ 1, e^{-t}, e^{-2t}, \dots, e^{-mt} \} \, , \\
U_\infty = \Ex (\cdot) \, ;
\end{gathered}
\end{equation}
the latter means that $ \lim_{t\to\infty} U_t X = \Ex(X) \cdot \One
$. According to \ref{1c}, $ U_t^A : L_2(\F_A) \to L_2(\F_A) $ for $ A
\in \A $ may be defined naturally, giving
\begin{equation}\label{1e}
\begin{gathered}
U_t^{A\uplus B} = U_t^A \otimes U_t^B \, , \qquad U_t^\emptyset =
  \One, \quad U_t^T = U_t \, , \\
U_\infty^A \otimes U_0^{T\setminus A} = \cE{ \cdot }{ \F_{T\setminus
  A} } \, , \\
\Proj_{H_A} = \bigg( \bigotimes_{a\subset A} \( \One - U_\infty^a \)
  \bigg) \otimes \bigg( \underbrace{ \bigotimes_{a\subset T\setminus
  A} U_\infty^a }_{=U_\infty^{T\setminus A}} \bigg) \, .
\end{gathered}\end{equation}
Introduce generators:
\begin{equation}\label{1f}
\begin{gathered}
U_t = \exp (-t\bN) \, , \quad \spec (\bN) = \{ 0,1,2,\dots,m \} \, , \\
X \in H_n \imply \bN X = n X \, ; \\
U_t^A = \exp (-t\bN^A) \, , \quad \bN^A : L_2 (\F_A) \to L_2 (\F_A) \,
, \\
\bN^T = \bN \, ; \qquad \bN_{a_k} X = X - \Ex X \quad
  \text{for } X \in L_2 (\F_{a_k}) \, ; \\
\bN^{A\uplus B} = \bN^A \otimes \One + \One \otimes \bN^B \, .
\end{gathered}
\end{equation}

The probabilistic meaning of $ U_t $ may be described roughly by
saying that each of our $ m $ pieces of data is unreliable, it is
either correct (with probability $ p $) or totally wrong (with
probability $ 1-p $). More exactly: any random variable $ X
\in L_2 (\F_T) $ is a function, $ X = \phi(Y_1,\dots,Y_m) $, of $ m $
random variables $ Y_1,\dots,Y_m $ such that $ Y_k $ is $ \F_{a_k}
$-measurable (therefore $ Y_1,\dots,Y_m $ are independent). Introduce
independent copies $ Z_1,\dots,Z_m $ of $ Y_1,\dots,Y_m $, and a
random set $ A \in \A $ distributed $ \mu_p $ and independent of $
Y_1,\dots,Y_m $, $ Z_1,\dots,Z_m $. Define $ Y'_1,\dots,Y'_m $ as
follows: if $ a_k \subset A $ then $ Y'_k = Y_k $, otherwise $ Y'_k =
Z_k $. We have
\begin{equation}\label{1.8}
\begin{gathered}
\cE{ \phi(Y_1,\dots,Y_m) }{ Y'_1,\dots,Y'_m } = \psi (Y'_1,\dots,Y'_m)
  \, , \\
U_t \( \phi (Y_1,\dots,Y_m) \) = \psi (Y_1,\dots,Y_m) \, ,
\end{gathered}
\end{equation}
which follows by averaging in $ A $ of
\begin{gather*}
\cE{ \phi(Y_1,\dots,Y_m) }{ A; \, Y'_1,\dots,Y'_m } = \psi_A
  (Y'_1,\dots,Y'_m) \, , \\
E_A \( \phi (Y_1,\dots,Y_m) \) = \psi_A (Y_1,\dots,Y_m) \, .
\end{gather*}
The reader may also imagine the corresponding continuous-time Markov
process; when $ a_k $ is excluded from our random set, the $ k $-th
portion of data is immediately replaced with an independent copy.
Such functions as\footnote{%
For their interrelations see the proof of Lemma \ref{2e}.}
$ t \mapsto \| X - U_t X \| $, $ t \mapsto \| X \| - \| U_t X \| $, or
$ t \mapsto \( (\One-U_t) X, X \) $ may be used for describing noise
sensitivity of a random variable $ X $. The more the functions, the
more sensitive is $ X $.
Least sensitive (most stable) are elements of $ H_1 = H_{a_1} \oplus
\dots \oplus H_{a_m} $, that is, random variables of the form $ X =
X_1 + \dots + X_m $, where $ X_k \in H_{a_k} $ (which means $ X_k \in
L_2 (\F_{a_k}) $, $ \Ex X_k = 0 $); these satisfy $ U_t X = e^{-t} X
$. Most sensitive are elements of $ H_m = H_T $, that is, linear
combinations of random variables of the form $ X = X_1 \dots X_m $
($ X_k $ being as above); these satisfy $ U_t X = e^{-mt} X $.  The
concept of noise sensitivity, quantitative for finite $ \A $, becomes
qualitative for infinite $ \A $, as we'll see in the next section.

The following result shows that contractions do not increase
sensitivity.
\begin{lemma}\label{1.10}
Let $ f : \R \to \R $ satisfy $ |f(x)-f(y)| \le |x-y| $ for all $ x,y
$. Then\footnote{%
For a random variable $ X : \Om \to \R $, $ f(X) $ denotes the
composition $ f \circ X : \Om \xrightarrow{X} \R \xrightarrow{f} \R
$.}
\[
\( (\One-U_t) f(X), f(X) \) \le \( (\One-U_t) X, X \)
\]
for all $ t \in [0,\infty) $ and $ X \in L_2(\F_T) $.
\end{lemma}

\begin{proof}
Introduce $ Y_1, \dots, Y_m $, $ Y'_1, \dots, Y'_m $, and $ \phi,\psi
$ as in \eqref{1.8}; note that $ Y'_1, \dots, Y'_m $ are independent
and distributed like $ Y_1, \dots, Y_m $; we have
\begin{multline*}
\Ex \( \phi (Y_1,\dots,Y_m) - \phi (Y'_1,\dots,Y'_m) \)^2 = \\
\Ex \( \phi (Y_1,\dots,Y_m) \)^2 + \Ex \( \phi (Y'_1,\dots,Y'_m) \)^2
  - 2 \Ex \( \phi (Y_1,\dots,Y_m) \phi (Y'_1,\dots,Y'_m) \) = \\
= \| X \|^2 + \| X \|^2 - 2 \Ex \( \phi (Y'_1,\dots,Y'_m) \cE{ \phi
  (Y_1,\dots,Y_m) }{ Y'_1,\dots,Y'_m } \) = \\
= 2 \| X \|^2 - 2 \Ex \( \phi (Y'_1,\dots,Y'_m) \psi (Y'_1,\dots,Y'_m)
  \) = \\
= 2 \| X \|^2 - 2 \( U_t X, X \) = 2 \( (\One-U_t) X, X \) \, ,
\end{multline*}
as well as
\[
\Ex \( f(\phi (Y_1,\dots,Y_m)) - f(\phi (Y'_1,\dots,Y'_m)) \)^2 = 2 \(
(\One-U_t) f(X), f(X) \) \, .
\]
However,\\
$ | f(\phi (Y_1,\dots,Y_m)) - f(\phi (Y'_1,\dots,Y'_m)) | \le
| \phi (Y_1,\dots,Y_m) - \phi (Y'_1,\dots,Y'_m) | $.
\end{proof}

Similarly, if $ f : \R^2 \to \R $ satisfies $ | f(x_1,x_2) -
f(y_1,y_2) | \le \( (x_1-y_1)^2 + (x_2-y_2)^2 \)^{1/2} $ then
\begin{equation}
\( (\One-U_t) f(X,Y), \, f(X,Y) \) \le \( (\One-U_t) X, \, X \) + \(
(\One-U_t) Y, \, Y \)
\end{equation}
for all $ X,Y \in L_2(\F_T) $ and $ t \in [0,\infty) $. The same for $
f : \R^d \to \R $.

\section{Stability, sensitivity, linearizability}
In this section the algebra $ \A $ is assumed to be countable. For
example, it may be the algebra generated by intervals $ (r,s) \subset
\R $ with rational $ r, s $, or the algebra of all cylindrical subsets
of $ \{0,1\}^\Z $. Being countable, $ \A $ is the union of a sequence
of its finite subalgebras:
\begin{equation}\label{2a}
\A = \A_1 \cup \A_2 \cup \dots \, , \qquad \A_1 \subset \A_2 \subset
\dots \text{ are finite subalgebras of $ \A $.}
\end{equation}
The freedom in choosing the sequence $ (\A_m) $ is of no importance
for us due to the following ``cofinality argument''. Let $ \phi $ be a
function defined on the set of all finite subalgebras of $ \A $ and
such that $ \lim_m \phi(\A_m) $ exists for every sequence $ (\A_m) $
satisfying \eqref{2a}. Then the limit is the same for all such
sequences. Proof: if $ (\A_m) $ and $ (\A'_m) $ are two such
sequences, then we can choose $ m_1 < m_2 < \dots $ and $ m'_1 < m'_2
< \dots $ such that
\begin{equation}\label{2b}
\A_{m_1} \subset \A'_{m'_1} \subset \A_{m_2} \subset \A'_{m'_2}
\subset \dots \, ,
\end{equation}
therefore the sequence $ \phi (\A_{m_1}), \phi (\A'_{m'_1}), \phi
(\A_{m_2}), \phi (\A'_{m'_2}), \dots $ must have a limit.

As before, $ (\F_A)_{A\in\A} $ satisfying \eqref{0.1} is
considered. (No other assumptions, such as \eqref{0.2}.) Still,
conditional expectation operators $ E_A $ are defined, see
\eqref{1.3}, \eqref{1.4}.

Restricting $ (\F_A) $ to $ A \in \A_m $ we get the elementary case
of Sect.~1. Probability measures $ \mu_p^{(m)} $ are defined, see
\eqref{1.1} and \eqref{1.2}, subspaces $ H_A^{(m)} $ for $ A \in \A_m
$, see \eqref{1.5},\footnote{%
The reader may guess that the decomposition of $ L_2 (\F_T) $ into the
direct sum of $ H_A^{(m)} $, $ A \in \A_m $, has a kind of limit for $
m \to \infty $. That is true; the limit is described in
\cite[Sect.~2]{unitary} in terms of direct integrals of Hilbert spaces
(for somewhat more restrictive framework, though). In the present
paper, direct integrals do not appear explicitly; however, most of the
text is in fact translated from that language.}
operator semigroups $ U_t^{(m)} $,
see \eqref{1.6}, their eigenspaces $ H_n^{(m)} $, see \eqref{1d}, and
generators $ \bN_m $, see \eqref{1f}. All $ U_t^{(m)} $ belong to the
commutative algebra generated by operators of conditional expectation
$ E_A $, $ A \in \A $. Compare $ H_A^{(m)} $ and $ H_B^{(m+1)} $ for $
A \in \A_m $, $ B \in \A_{m+1} $; the second space is either included
into the first, or orthogonal to it; namely, if $ A $ is the least
element of $ \A_m $ containing $ B $ ($ \A_m $-saturation of $ B $),
then $ H_B^{(m+1)} \subset H_A^{(m)} $, otherwise $ H_B^{(m+1)} \bot
H_A^{(m)} $. If $ X \in H_B^{(m+1)} \subset H_A^{(m)} $ then $
U_t^{(m)} X = e^{-kt} X $, $ U_t^{(m+1)} X = e^{-lt} X $ with $ k \le
l $ (since saturation does not increase the number of atoms). So,
\begin{equation}\label{2c}
\bN_m \le \bN_{m+1} \, ; \qquad U_t^{(m)} \ge U_t^{(m+1)} \, .
\end{equation}
It follows easily that the limit exists,
\begin{equation}\label{2d}
\begin{gathered}
U_t = \lim_{m\to\infty} U_t^{(m)} \quad \text{in the sense that }
  \forall X \in L_2 \;\; \| U_t^{(m)} X - U_t X \| \to 0 \, ; \\
U_s U_t = U_{s+t} \, ; \qquad \| U_t \| \le 1 \, .
\end{gathered}\end{equation}
The limit, $ U_t $, does not depend on the choice of $ (\A_m) $ due to
the cofinality argument (see \eqref{2b}). Also, $ U_t $ commute with
all $ E_A $.

The limit of generators, $ \lim_m \bN_m $, need not exist; $ \| \bN_m
X \| $ can tend to $ \infty $ for some $ X $. Accordingly, the
operator semigroup $ (U_t) $ need not be continuous at $ t = 0 $.

\begin{lemma}\label{2e}
There exists a sub-\sif\ $ \Fstable \subset \F_T $ such that
\[
X \in L_2 (\Fstable) \equiv \| X - U_t X \| \xrightarrow[t\to0]{} 0 \, .
\]
\end{lemma}

\begin{proof}
First, the following three properties of $ X $ are equivalent: (a) $
\| X - U_t X \| \xrightarrow[t\to0]{} 0 $; (b) $ \| U_t X \|
\xrightarrow[t\to0]{} \| X \| $; (c) $ \( (\One-U_t) X, X \)
\xrightarrow[t\to0]{} 0 $. Indeed, \mbox{(c) \imp (a)} since $ \| X -
U_t X \|^2 = \( (\One-U_t)^2 X, X \) \le \( (\One-U_t) X, X \) $;
\mbox{(a) \imp (b)} since $ \| U_t X \| \ge \| X \| - \| X - U_t X \|
$; \mbox{(b) \imp (c)} since $ \( (\One-U_t) X, X \) = \| X \|^2 - \|
U_{t/2} X \|^2 $.

The set $ H_{\text{stable}} = \{ X \in L_2 (\F_T) : \| X - U_t X \|
\xrightarrow[t\to0]{} 0 \} = \{ X \in L_2 (\F_T) : \( (\One-U_t) X, X \)
\xrightarrow[t\to0]{} 0 \} $ is a closed linear subspace of $ L_2 $. By
Lemma \ref{1.10}, if $ X,Y \in H_{\text{stable}} $ then $ \min(X,Y),
\max(X,Y) \in H_{\text{stable}} $. Also, $ H_{\text{stable}} $
contains constants. It is well-known that such a space is the whole $
L_2 (\Fstable) $ where $ \Fstable $ is the \sif\ generated by $
H_{\text{stable}} $.
\end{proof}

We have
\begin{equation}\label{2f}
\begin{split}
U_t \( L_2 (\Fstable) \) &\subset L_2 (\Fstable) \quad \text{for all }
  t \in [0,\infty) \, , \\
E_A \( L_2 (\Fstable) \) &\subset L_2 (\Fstable) \quad \text{for all }
  A \in \A \, ,
\end{split}\end{equation}
since $ \| U_s X - U_t U_s X \| = \| U_s ( X - U_t X ) \| \le \| X -
U_t X \| $ and $ \| E_A X - U_t E_A X \| = \| E_A ( X - U_t X ) \| \le
\| X - U_t X \| $.
Being restricted to $ L_2 (\Fstable) $, the operator semigroup $ (U_t)
$ is continuous (in the strong operator topology) and has its
generator $ \bN = \lim_m \bN_m $, $ \spec \bN \subset \{ 0,1,2,\dots
\} $; denote its eigenspaces by $ H_n $;
\begin{equation}\label{2g}
\begin{gathered}
U_t = e^{-t\bN} \quad \text{on } L_2 (\Fstable) \, , \\
L_2 (\Fstable) = H_0 \oplus H_1 \oplus H_2 \oplus \dots \, , \\
U_t = e^{-nt} \quad \text{on } H_n \, , \\
E_A (H_n) \subset H_n \quad \text{for all } A \in \A \, ;
\end{gathered}\end{equation}
check the latter: $ X \in H_n \impl U_t E_A X = E_A U_t X = e^{-nt}
E_A X \impl E_A X \in H_n $.
The relation $ \bN_m \uparrow \bN $ implies for all $ n \in
\{0,1,2,\dots\} $
\begin{equation}\label{2h}
H_0 \oplus \dots \oplus H_n = \bigcap_{m=1}^\infty H_0^{(m)} \oplus
\dots \oplus H_n^{(m)} 
\end{equation}
(the intersection of a decreasing sequence of subspaces). Clearly, $
H_0 = H_0^{(m)} = H_\emptyset^{(m)} $ is the one-dimensional space of
constants.

\begin{lemma}\label{2i}
The following conditions are equivalent for all $ X \in L_2 (\F_T)
$:

(a) $ X \in H_1 $;

(b) $ X = E_A X + E_{T\setminus A} X $ for all $ A \in \A $;

(c) $ X = E_{A_1} X + \dots + E_{A_k} X $ for every partition $ T = A_1
\uplus \dots \uplus A_k $ of $ T $ into $ A_i \in \A $.
\end{lemma}

\begin{proof}
Each element of $ H_1^{(m)} $ satisfies (c)
for the partition into atoms of $ \A_m $. Therefore each element of $
H_1 = \cap_m H_1^{(m)} $ satisfies (c) for every partition, and we get
(a) \imp (c) \imp (b). For proving (b) \imp (a) assume that $ X $
satisfies (b) and prove that $ X \in H_1^{(m)} $ for all $ m $. From
now on $ A $ and $ B $ run over $ \A_m $. By \eqref{1.5}, $ X = \sum
X_A $, $ X_A \in H_A^{(m)} $; we have to prove that $ X_A = 0 $ unless
$ A $ contains exactly one atom. By \eqref{1.5} again,
\[
E_A X = \sum_{B\subset A} X_B \, , \qquad E_{T\setminus A} X =
\sum_{B\subset T\setminus A} X_B \, . 
\]
We see that $ X_\emptyset $ appears twice in $ E_A X + E_{T\setminus
A} X $, but only once in $ X $, therefore $ X_\emptyset = 0 $. If $
B $ contains at least two atoms, we can choose $ A $ such that $ B $
intersects both $ A $ and $ T \setminus A $; then $ X_B $ does not
appear in $ E_A X + E_{T\setminus A} X $, but appears in $ X $,
therefore $ X_B = 0 $.
\end{proof}

The following is a general fact about Hilbert spaces, irrespective of
any probability theory.

\begin{lemma}\label{2j}
Assume that $ H' $ and $ H'' $ are Hilbert spaces, $ H = H' \otimes
H'' $, and subspaces are given, $ H' \supset H'_1 \supset H'_2 \supset
\dots $ and $ H'' \supset H''_1 \supset H''_2 \supset \dots \, $ Then
\[
\bigcap_m ( H'_m \otimes H''_m ) = \bigg( \bigcap_m H'_m \bigg)
\otimes \bigg( \bigcap_m H''_m \bigg) \, .
\]
\end{lemma}

\begin{proof}
Denoting $ H'_0 = H' $, $ H'_\infty = \cap H'_m $, we have
\[
H' = H'_\infty \oplus \bigoplus_{m=0}^\infty ( H'_m \ominus H'_{m+1} )
\]
and the same for $ H'' $. Therefore
\begin{multline*}
H = \Big( H'_\infty \otimes H''_\infty \Big) \oplus
  \bigg( \bigoplus_{m=0}^\infty ( H'_m
  \ominus H'_{m+1} ) \otimes H''_\infty \bigg) \oplus \\
\oplus \bigg( H'_\infty \otimes
  \bigoplus_{n=0}^\infty ( H''_n \ominus H''_{n+1} ) \bigg) \oplus
  \bigg( \bigoplus_{m,n} ( H'_m \ominus H'_{m+1} ) \otimes ( H''_n
  \ominus H''_{n+1} ) \bigg) \, .
\end{multline*}
The space $ H'_m \otimes H''_m $ contains some of the terms, and is
orthogonal to others. Only the term $ H'_\infty \otimes H''_\infty $
is contained in $ H'_m \otimes H''_m $ for all $ m $.
\end{proof}

\begin{lemma}\label{2k}
For any $ m $ and any two different atoms $ a,b $ of $ \A_m $,
\[
H_2 \cap H_{a\cup b}^{(m)} = \( H_1 \cap H_a^{(m)} \) \otimes \( H_1
\cap H_b^{(m)} \) \, .
\]
\end{lemma}

\begin{proof}
$ H_{a\cup b}^{(m)} = H_a^{(m)} \otimes H_b^{(m)} $ by \eqref{1a}; $
H_2 \cap H_{a\cup b}^{(m)} = ( H_0 \oplus H_1 \oplus H_2 ) \cap
H_{a\cup b}^{(m)} = \cap_k \( H_0^{(m+k)} \oplus H_1^{(m+k)} \oplus
H_2^{(m+k)} \) \cap H_{a\cup b}^{(m)} = \cap_k H_2^{(m+k)} \cap
H_{a\cup b}^{(m)} $ by \eqref{2h}; note that the space decreases when
$ k $ increases. Similarly, $ H_1 \cap H_a^{(m)} = \cap_k H_1^{(m+k)}
\cap H_a^{(m)} $ and $ H_1 \cap H_b^{(m)} = \cap_k H_1^{(m+k)} \cap
H_b^{(m)} $. By Lemma \ref{2j} it suffices to prove that
\[
H_2^{(m+k)} \cap H_{a\cup b}^{(m)} = \( H_1^{(m+k)} \cap H_a^{(m)} \)
\otimes \( H_1^{(m+k)} \cap H_b^{(m)} \)
\]
for $ k = 1,2,\dots $
However, $ H_2^{(m+k)} $ is (by definition) the direct sum of $
H_{c\cup d}^{(m+k)} $ over atoms $ c,d $ of $ \A_{m+k} $, $ c \ne d $,
and $ H_2^{(m+k)} \cap H_{a\cup b}^{(m)} $ is such a sum over $ c
\subset a $, $ d \subset b $. It remains to note that $ H_{c\cup
d}^{(m+k)} = H_c^{(m+k)} \otimes H_d^{(m+k)} $.
\end{proof}

\begin{theorem}\label{2l}
The \sif\ generated by $ H_1 $ is equal to $ \Fstable $.
\end{theorem}

\begin{proof}
Denote by $ \F_n $ the \sif\ generated by $ H_n $. It suffices to
prove that $ \F_n \subset \F_1 $ for all $ n $, since $ L_2 (\Fstable)
= H_0 \oplus H_1 \oplus H_2 \oplus \dots $ (see \eqref{2g}). I give a
proof for $ n = 2 $; it has a straightforward generalization for
higher $ n $.

We have to prove that $ H_2 \subset L_2 (\F_1) $. For each $ m $, $
H_2 = \( H_2 \cap H_1^{(m)} \) \oplus \( H_2 \cap H_2^{(m)} \) $
(since $ H_2 $ is
invariant under all $ E_A $). However, $ H_1^{(m)} $ decreases to $
H_1 $, and $ H_1 $ is orthogonal to $ H_2 $. Therefore the union of $
H_2 \cap H_2^{(m)} $ is dense in $ H_2 $; it remains to prove that $
H_2 \cap H_2^{(m)} \subset L_2 (\F_1) $ for all $ m $. Note that $ H_2
\cap H_2^{(m)} $ is the direct sum of $ H_2 \cap H_{a\cup b}^{(m)} $
over atoms $ a,b $ of $ \A_m $, $ a \ne b $ (since, again, $ H_2 $ is
invariant under all $ E_A $). Lemma \ref{2k} reduces the needed
inclusion to an evident fact, $ ( H_1 \cap H_a^{(m)} ) \otimes ( H_1
\cap H_b^{(m)} ) \subset L_2 (\F_1) $.
\end{proof}

A canonical isomorphism between $ H_n $ and $ \underbrace{ H_1 \otimes
  \dots \otimes H_1 }_n $ is given by Wick products,\footnote{%
Given $ X,Y \in H_1 $, we may define their Wick product,
\[
{:}XY{:} \, = \, \lim_{m\to\infty} \sum_{a\ne b} (E_a X) (E_b Y) \, .
\]
The sum is taken over all unordered pairs $ \{ a,b \} $ of different
atoms of $ \A_m $. The same for $ {:}XYZ{:} $ and so on.}
but is not needed here.

\begin{definition}\label{2m}
(a) A random variable $ X \in L_2 (\F_T) $ is called \emph{stable,} if
$ X \in L_2 (\Fstable) $, and \emph{sensitive,} if $ \cE{ X }{
\Fstable } = 0 $.

(b) Equivalently, a random variable $ X \in L_2 (\F_T) $ is called
\emph{stable,} if $ \| X - U_t X \| \xrightarrow[t\to0]{} 0 $, and
\emph{sensitive,} if $ U_t X = 0 $ for all $ t \in (0,\infty) $.
\end{definition}

The two definitions of stability are equivalent evidently (recall
\ref{2e}), of sensitivity --- due to the following result.

\begin{lemma}\label{2n}
The following conditions are equivalent for all $ X \in L_2 (\F_T) $:

(a) $ \cE{ X }{ \Fstable } = 0 $;

(b) $ U_t X = 0 $ for all $ t > 0 $.
\end{lemma}

\begin{proof}
(b) \imp (a):
Let $ Y \in L_2 (\Fstable) $, then $ (X,Y) = \lim_{t\to0} (X,U_t Y) =
\lim_{t\to0} (U_t X, Y) = 0 $ by (b).

(a) \imp (b):
It suffices to prove that $ \| U_t X \| \le e^{-nt} \| X \| $ for all
$ X \in L_2 (\F_T) $ orthogonal to $ H_0 \oplus \dots \oplus H_{n-1}
$. By \eqref{2h} we may assume that $ X $ is orthogonal to $ H_0^{(m)}
\oplus \dots \oplus H_{n-1}^{(m)} $ for some $ m $ (since such vectors
are dense in $ L_2 (\F_T) \ominus ( H_0 \oplus \dots \oplus H_{n-1} )
$). For such $ X $, $ \| U_t X \| \le \| U_t^{(m)} X \| \le e^{-nt} \|
X \| $.
\end{proof}

So, in terms of $ U_{0+} X = \lim_{t\to0,t>0} U_t X $ we have
\begin{equation}\label{2o}
\begin{gathered}
L_2 (\F_T) = \{ X : X \text{ is stable} \, \} \oplus \{ X : X \text{ is
  sensitive} \, \} \, , \\
X \text{ is stable } \equiv U_{0+} X = X \, , \\
X \text{ is sensitive } \equiv U_{0+} X = 0 \, ; \\
\cE{ \cdot }{ \Fstable } = U_{0+} \, .
\end{gathered}
\end{equation}
Similarly, for any $ A \in \A $
\begin{equation}\label{2p}
\cE{ \cdot }{ \Fstable^A } = U_{0+}^A \quad \text{on } L_2 (\F_A)
\end{equation}
for some $ \Fstable^A $, and $ \Fstable^A = \F_A \cap \Fstable $,
since $ U_t^A X = U_t X $ for $ X \in L_2 (\F_A) $; also, $ \Fstable^A
$ is generated by $ H_1^A = H_1 \cap L_2 (\F_A) $, therefore
\begin{equation}\label{2q}
\Fstable^{A\uplus B} = \Fstable^A \otimes \Fstable^B \, ,
\end{equation}
which means that $ \( \Fstable^A \)_{A\in\A} $ is another family
satisfying \eqref{0.1}, the \emph{stable (or linearizable) part} of
the given family $ (\F_A)_{A\in\A} $. (See also \cite[Th.~1.7]{TV}.)

\section{Stability and extendibility}
We still work with a countable algebra $ \A $ and a family $
(\F_A)_{A\in\A} $ satisfying \eqref{0.1}. Striving to extend the
family from the algebra $ \A $ to the \sif\ generated by $ \A $ we can
face the following obstacle.

Let $ A_k \in \A $, $ A_1 \subset A_2 \subset \dots $; consider two
\sif s: $ \bigvee \F_{A_k} $ (the least \sif\ containing all $
\F_{A_k} $), and $ \bigwedge \F_{T\setminus A_k} = \bigcap
\F_{T\setminus A_k} $ (the intersection of all $ \F_{T\setminus A_k}
$). It is easy to see that the two \sif s are independent. The
question is, whether
\begin{equation}\label{3a}
\bigg( \bigvee_k \F_{A_k} \bigg) \otimes \bigg( \bigwedge_k
\F_{T\setminus A_k} \bigg) = \F_T \, ,
\end{equation}
or not. That is, whether the two \sif s generate the whole $ \F_T $,
or not. If they do not, then $ (\F_A) $ has no $\sigma$-additive (in
the sense of \eqref{0.2}) extension to a \sif.

\begin{theorem}\label{3b}
If \eqref{3a} is satisfied for every increasing sequence $ (A_k) $,
then $ \F_T = \Fstable $.
\end{theorem}

Postpone the proof. Choose $ \A_m $ satisfying \eqref{2a}. Choose $
p_k \in (0,1) $ such that $ \sum (1-p_k) < 1 $, say, $ p_k = 1 -
2^{-k-1} $. Recall probability measures $ \mu_p^{(m)} $.

\begin{lemma}\label{3c}
There exists a sequence $ m_1 < m_2 < \dots $ such that $ \(
\mu_{p_1}^{(m_1)} \otimes \mu_{p_2}^{(m_2)} \otimes \dots \) $-almost
all sequences $ (A_1,A_2,\dots) $, $ A_k \in \A_{m_k} $, satisfy
\[
\bigcap_{k=1}^\infty \F_{A_k} \subset \Fstable \, .
\]
\end{lemma}

\begin{proof}
If $ \bigcap \F_{A_k} $ is not contained in $ \Fstable $ then there
exists $ X \in L_2 ( \bigcap \F_{A_k} ) $, $ X \ne 0 $, orthogonal
to $ L_2 (\Fstable) $, that is, sensitive. The case is impossible,
if $ E_{A_k} X \to 0 $ for all sensitive $ X $ or, equivalently, for
a dense set of such $ X $; the more so, if $ \sum_k (E_{A_k} X, X) <
\infty $ for all these $ X $. By \eqref{1.6}, $ ( U_{t_k}^{(m_k)} X,
X ) $ is the average of $ ( E_{A_k} X, X ) $ over $ A_k $
distributed $ \mu_{p_k}^{(m_k)} $; here $ t_k = - \ln p_k $. It
suffices to choose $ m_k $ such that $ \sum_k ( U_{t_k}^{(m_k)} X, X
) < \infty $ for a dense set of sensitive $ X $.

For each sensitive $ X $ and each $ t > 0 $, by \ref{2n}, $ (
U_t^{(m)} X, X ) \to 0 $ for $ m \to \infty $. Therefore $
(U_{t_k}^{(m_k)} X, X) \to 0 $ for $ k \to \infty $, if $ m_k $ grow
fast enough. Diagonal argument gives a single sequence $ (m_k) $ that
serves a given sequence of vectors $ X $. It remains to choose a
sequence dense among all sensitive vectors.
\end{proof}

Introduce $ \tilde E_A $ similar to $ E_A $ as follows:
\begin{equation}\label{3d}
\begin{split}
E_A &= \cE{ \cdot }{ \F_A } = \One_A \otimes U_\infty^{T\setminus A}
  \, ; \\
\tilde E_A &= \cE{ \cdot }{ \F_A \vee \Fstable } = \One_A \otimes
  U_{0+}^{T\setminus A} \, ;
\end{split}\end{equation}
here $ \F_A \vee \Fstable $ is the \sif\ generated by these two \sif
s, and $ \One_A $ is the unit operator on $ L_2 (\F_A) $; the equality
$ \cE{ \cdot }{ \F_A \vee \Fstable } = \One_A \otimes
U_{0+}^{T\setminus A} $ follows from \eqref{2p}, since by
\eqref{2q}, $ \F_A \vee \Fstable = \F_A \vee \( \Fstable^A \otimes
\Fstable^{T\setminus A} \) = \F_A \otimes \Fstable^{T\setminus A} $.

We could proceed to $ \tilde U_t^{(m)} $ similar to $ U_t^{(m)} $,
\begin{equation}\label{3e}
\begin{split}
U_t^{(m)} = \int E_A \, d\mu_p^{(m)} (A) \, , \\
\tilde U_t^{(m)} = \int \tilde E_A \, d\mu_p^{(m)} (A) \, ,
\end{split}
\qquad
(p=e^{-t})
\end{equation}
and to $ \tilde U_t = \lim_m U_t^{(m)} $; however, we need a bit more
general construction,
\begin{equation}\label{3f}
\tilde U_\mu = \int \tilde E_A \, d\mu (A)
\end{equation}
for an arbitrary probability distribution $ \mu $ on $ A $ ($ A $ is
treated here as just a countable set). In fact, we need only $ \mu $
concentrated on a finite set, which is elementary in the sense of
Sect.~1.

\begin{lemma}\label{3g}
$ \tilde U_\mu \le (1-p) U_{0+} + p \cdot \One $, where\footnote{%
Recall that $ \A $ is an algebra of subsets of some set $ T $. The
latter was mentioned only once, before \eqref{0.1}, and may be readily
avoided now;
\[
p = \sup_{B\in\A,B\ne\emptyset} \mu \( \{ A \in \A : A \supset B \} \)
\, .
\]
}
\[
p = \sup_{t\in T} \mu \( \{ A \in \A : A \ni t \} \) \, .
\]
\end{lemma}

\begin{proof}
Similarly to \eqref{1b1}, \eqref{1e}, for every $ m $,
\begin{equation}\label{3h}
\begin{gathered}
L_2 (\F_T) = \bigotimes_a L_2 (\F_a) = \bigotimes_a \(
  H_{\text{stable}}^a \oplus H_{\text{sensitive}}^a \) =
  \bigoplus_{A\in\A_m} \tilde H_A \, , \\
\Proj_{\tilde H_A} = \bigg( \bigotimes_{a\subset A}
  \underbrace{ \Proj_{H_{\text{sensitive}}^a} }_{\One - U_{0+}^a}
 \bigg) \otimes \bigg(
  \bigotimes_{a\subset T\setminus A} \underbrace{
  \Proj_{H_{\text{stable}}^a} }_{U_{0+}^a} \bigg) \, ,
\end{gathered}
\end{equation}
where $ a $ runs over atoms of $ \A_m $, and $ H_{\text{stable}}^a $,
$ H_{\text{sensitive}}^a $ are subspaces of stable and sensitive,
respectively, elements of $ L_2 (\F_a) $.

For every $ B \in \A_m $, $ B \ne \emptyset $, the operator $ \tilde
E_A $ (recall \eqref{3d}) on $ H_B $ is the unit (identity) if $ B
\subset A $, otherwise it vanishes. Assuming $ \mu (\A_m) = 1 $ we get
$ \tilde U_\mu = \la \cdot \One $ on $ H_B $, where $ \la = \mu \( \{
A : A \supset B \} \) \le p $. Therefore $ \tilde U_\mu \le p \cdot
\One $ on $ H_B $, $ B \ne \emptyset $ (note that $ \tilde U_\mu (H_B)
\subset H_B $), while on $ H_\emptyset $ we have $ U_{0+} = \One $;
so, $ \tilde U_\mu \le (1-p) U_{0+} + p \cdot \One $ provided that $
\mu (\A_m) = 1 $ for some $ m $. The general case, $ \mu (\A_m) \to 1
$, will not be used, and I leave it to the reader.
\end{proof}

\begin{proof}[Proof of Theorem \ref{3b}]
Choose $ m_k $ by Lemma \ref{3c}, then $ \cap \F_{A_k} \subset
\Fstable $ for $ \mu $-almost all $ (A_k) $; here $ \mu = \otimes_k
\mu_{p_k}^{(m_k)} $. On the other hand, for every $ t \in T $ and
every $ k $,
\begin{multline*}
\mu \( \{ (A_k) : t \in T \setminus (A_1 \cap\dots\cap A_k) \} \) \le
  \\
\le \mu_{p_1}^{(m_1)} \( \{ A_1 : t \in T \setminus A_1 \} \) + \dots
  + \mu_{p_k}^{(m_k)} \( \{ A_k : t \in T \setminus A_k \} \) \le \\
  \le \sum_i (1-p_i) = q < 1 \, ;
\end{multline*}
by Lemma \ref{3g},
\[
\int \tilde E_{T\setminus(A_1\cap\dots\cap A_k)} \, d\mu \le (1-q)
  U_{0+} + q \cdot \One \, ,
\]
therefore
\begin{multline*}
\int \bigg\| \CE{ X }{ \bigvee_{k=1}^\infty \F_{T \setminus (A_1
  \cap\dots\cap A_k)} \vee \Fstable } \bigg\|^2 \, d\mu = \\
= \lim_{k\to\infty} \int \| \cE{ X }{ \F_{T\setminus (A_1
  \cap\dots\cap A_k)} \vee \Fstable } \|^2 \, d\mu \le q \| X \|^2
\end{multline*}
for all sensitive $ X \in L_2 (\F_T) $. Applying \eqref{3a} to the
increasing sequence $ T \setminus (A_1 \cap\dots\cap A_k) $ we get  
\[
\bigg( \bigvee_{k=1}^\infty \F_{T\setminus (A_1 \cap\dots\cap A_k)}
\bigg) \otimes \bigg( \bigwedge_{k=1}^\infty \F_{A_k} \bigg) = \F_T
\]
for all $ (A_k) $, therefore
\[
\bigg( \bigvee_{k=1}^\infty \F_{T\setminus (A_1 \cap\dots\cap A_k)}
\bigg) \vee \Fstable = \F_T
\]
for $ \mu $-almost all $ (A_k) $. So, each sensitive $ X $ satisfies $
\int \| X \|^2 \, d\mu \le q \| X \|^2 $, that is, $ \| X \|^2 \le q
\| X \|^2 $, which is impossible unless $ X = 0 $.
\end{proof}

So, if $ \F_T \ne \Fstable $ then $ (\F_A)_{A\in\A} $ has no
$\s$-additive extension. On the other hand, if $ \F_T = \Fstable $
then such an extension is usually possible, for a simple reason: $ E_A
$ restricted to $ H_1 $ form a projection-valued finitely additive
measure. Conditions well-known in measure theory ensure that a
$\s$-additive extension to a \sif\ exists, and we get extended $ \F_A
$ as generated by extended $ H_1^A $.

\appendix
\section{Appendix: The simplest example of sensitivity}
\vspace{-0.5cm}
\hfill\parbox{7cm}{%
The phenomenon \dots tripped up even Kolmogorov
  and Wiener. \cite[p.~48]{Wi}
}
\vspace{0.5cm}

Two examples of a countable algebra $ \A $, mentioned in the beginning
of Sect.~2, are nonatomic; corresponding families $ (\F_A)_{A\in\A} $
are in general as complicated as continuous-time random processes. The
simplest infinite $ \A $ consists of all finite and
cofinite\footnote{%
A set is called cofinite if its complement is finite.}
subsets of $ T = \{ 1,2,\dots \} $. From now on, $ \A $ stands for
that algebra; it is purely atomic, and corresponding $ (\F_A)_{A\in\A}
$ are as simple as discrete-time random processes, that is, random
sequences. Not too simple, as we'll see soon\dots

Choose some $ p \in \{ 3,5,7,9,\dots \} $ and consider the simple
stationary random walk on the finite group $ \Z_p $. That is, $ \Om $
is the set of all sequences $ \om = (x_0,x_1,x_2,\dots) $, $ x_k \in
\Z_p $, $ x_{k+1}-x_k = \pm 1 $; $ \F $ is the \sif\ generated by
cylinder sets $ E_{y_0,\dots,y_m} = \{ \om \in \Om : X_0(\om) = y_0,
\dots, X_m(\om) = y_m \} $, where $ X_k (x_0,x_1,\dots) = x_k $; and $
P $ is defined by $ P (E_{y_0,\dots,y_m} = p^{-1}2^{-m} $ whenever $
y_k \in \Z_p $, $ y_{k+1}-y_k = \pm 1 $. So, each of the $ \Z_p
$-valued random variables $ X_0, X_1, \dots $ is uniformly
distributed; increments $ X_1 - X_0, X_2 - X_1, \dots $ are
independent, $ \pm 1 $ with probabilities $ 1/2, 1/2 $; and the random
variables $ X_0; X_1-X_0, X_2-X_1, \dots $ are independent.

Define \sif s $ \F_A $ for $ A \in \A $:
\begin{equation}\label{Aa}
\begin{split}
\F_{ \{k\} } &= \s (X_k-X_{k-1}) \, , \\
\F_{ \{k,k+1,\dots\} } &= \s ( X_{k-1}, X_k, X_{k+1}, \dots ) \, , \\
\F_{ \{k_1,\dots,k_n\} } &= \F_{k_1} \vee \dots \vee \F_{k_n} \, , \\
\F_{ \{k_1,\dots,k_n\} \cup \{k,k+1,\dots\} } &= \F_{ \{k_1,\dots,k_n\}
  } \vee \F_{ \{k,k+1,\dots\} } \, ;
\end{split}\end{equation}
here $ n \in \{0,1,\dots\} $, $ k, k_1, \dots, k_n \in \{1,2,\dots\}
$, $ k_1 < \dots < k_n < k $, and $ \s(\dots) $ means the \sif\
generated by given random variables. It is not immediately clear that
the definition is correct and \eqref{0.1} is satisfied, but it is
true; you may check it, starting with
\[
\F_{ \{1,\dots,k-1\} } \otimes \F_{ \{k,k+1,\dots\} } = \F \, .
\]

\begin{sloppypar}
Condition \eqref{3a} is violated for $ A_k = \{1,\dots,k\} $, since
the \sif\ $ \bigwedge_{k=1}^\infty \F_{ \{k+1,k+2,\dots\} } $ is
degenerate, while the \sif\ $ \bigvee_{k=1}^\infty \F_{ \{1,\dots,k\}
} = \s ( X_1-X_0, X_2-X_1, \dots ) $ contains only sets invariant
under the symmetry
\begin{equation}\label{Ab}
R : \Om \to \Om \, , \qquad R (x_0,x_1,\dots) = (x_0+1, x_1+1, \dots )
\, .
\end{equation}
(Note that $ X_0 $ is not invariant under $ R $.) Therefore $ \F_A $
cannot be defined for all $ A \subset \{1,2,\dots\} $ obeying
\eqref{0.1}, \eqref{0.2} and \eqref{Aa}.
\end{sloppypar}

We choose finite subalgebras $ \A_m \subset \A $, satisfying
\eqref{2a}, in a natural way:
\[
\text{atoms of $ \A_m $ are} \quad \{1\}, \, \dots, \, \{m-1\}, \,
\text{ and } \{m,m+1,\dots\} \, .
\]
An elementary calculation, starting with
\[
X_0 = X_{m-1} - (X_{m-1}-X_{m-2}) - \dots - (X_1-X_0) \, ,
\]
gives
\[\begin{split}
& U_t^{(m)} \exp \bigg( \frac{2\pi i}{p} X_0 \bigg) = \\
& \quad = e^{-t} \exp \bigg( \frac{2\pi i}{p} X_{m-1} \bigg) \cdot
  \prod_{k=1}^{m-1} \bigg( \cos \frac{2\pi}p + i e^{-t} \sin
  \frac{2\pi}p (X_k-X_{k-1}) \bigg) \, ; \\
& \bigg\| U_t^{(m)} \exp \bigg( \frac{2\pi i}{p} X_0 \bigg) \bigg\| =
  e^{-t} \bigg( \cos^2 \frac{2\pi}p + e^{-2t} \sin^2 \frac{2\pi}p
  \bigg)^{(m-1)/2} \, ;
\end{split}\]
therefore $ U_t \exp \( \frac{2\pi i}{p} X_0 ) = 0 $ for all $ t > 0
$, which means that
\begin{equation}
\exp \bigg( \frac{2\pi i}{p} X_0 \bigg) \quad \text{is sensitive.}
\end{equation}
In fact, $ \Fstable = \s (X_1-X_0, X_2-X_1, \dots) $ is the \sif\ of
all measurable sets that are invariant under the symmetry $ R
$. Accordingly,
\begin{equation}
\cE{ X }{ \Fstable } = \frac1p ( X + X \circ R + X \circ R^2 + \dots +
X \circ R^{p-1} ) \, .
\end{equation}
Also, it is easy to see that
\begin{equation}
H_1 = \{ c_1(X_1-X_0) + c_2(X_2-X_1) + \ldots \, : \, (c_1,c_2,\dots)
\in l_2
\}
\end{equation}
(here $ X_k - X_{k-1} $ is treated as taking on values $ \pm 1 \in \R
$ rather than $ \pm 1 \in \Z_p $).

Instead of $ \Z_p $ we could consider the unit circle on the complex
plane, and some random walk in the circle (or another compact group).

A physicists could write
\[
\exp \bigg( \frac{2\pi i}{p} X_0 \bigg) = \bigotimes_{k=1}^\infty \exp
\bigg( \frac{2\pi i}{p} (X_k-X_{k-1}) \bigg) 
\]
and say: that is just the wave function of an infinite sequence of
uncorrelated spins (or quantum bits), all in the same superposition of
two basis states. True, the infinite product of independent
identically distributed random variables does not converge, but
anyway, infinitely many commuting copies of $ \mathrm{SU(2)} $ act on
$ L_2 (\F_T) $.

\bigskip
\filbreak
\begingroup
{
\small
\begin{sc}
\parindent=0pt\baselineskip=12pt

School of Mathematics, Tel Aviv Univ., Tel Aviv
69978, Israel
\emailwww{tsirel@math.tau.ac.il}
{http://math.tau.ac.il/$\sim$tsirel/}
\end{sc}
}
\filbreak

\endgroup

\end{document}